\documentclass{amsproc}
\usepackage{graphicx,mathrsfs,amssymb,amsmath}
 \usepackage{xypic}

\newtheorem{theorem}{Theorem}[section]

\newtheorem{proposition}[theorem]{Proposition}
\newtheorem{remark}[theorem]{Remark}
\newtheorem{definition}[theorem]{Definition}

\theoremstyle{remark}

\numberwithin{equation}{section}

\newcommand{\cz}{{\mathbb C}}
\newcommand{\gz}{{\mathbb Z}}

\newcommand{\nz}{{\mathbb N}}

\newcommand{\rz}{{\mathbb R}}

\newcommand{\calA}{\mathcal{A}}

\newcommand{\calD}{\mathcal{D}}

\newcommand{\calU}{\mathcal{U}}

\newcommand{\scrC}{\mathscr{C}}
\newcommand{\scrD}{\mathscr{D}}

\newcommand{\PsDO}{\Psi \mathrm{DO}}

\newcommand{\dbar}{d\hspace*{-0.08em}\bar{}\hspace*{0.1em}}

\newcommand{\forget}[1]{}

\newcommand{\lra}{\longrightarrow}
\newcommand{\op}{\mathrm{op}}

\newcommand{\res}{\mathrm{res}}

\newcommand{\skp}[2]{\langle#1,#2\rangle}

\newcommand{\tr}{\mathrm{tr}\,}

\newcommand{\wt}{\widetilde}

\begin{document}
\title[Noncommutative residue for projective $\psi$do]{
On the Noncommutative Residue for \\ Projective Pseudodifferential Operators}

\author{J\"org Seiler}
\address{Loughborough University, Department of Mathematical Sciences,
         Leicestershire LE11 3TU (UK)}
\email{j.seiler@lboro.ac.uk}
\author{Alexander Strohmaier}
\address{Loughborough University, Department of Mathematical Sciences,
         Leicestershire LE11 3TU (UK)}
\email{a.strohmaier@lboro.ac.uk}

\begin{abstract}
A well known result on pseudodifferential operators states that the 
noncommutative residue (Wodzicki residue) of a pseudodifferential
projection vanishes. This statement is non-local and implies the
regularity of the eta invariant at zero of Dirac type operators.
We prove that in a filtered algebra the value of a projection under any residual
trace depends only on the principal part of the projection. This
general, purely algebraic statement applied to the algebra of projective pseudodifferential operators
implies that the noncommutative residue factors
to a map from the twisted $K$-theory of the co-sphere bundle. We use arguments
from twisted $K$-theory to show that this map vanishes in odd dimensions, thus showing that
the noncommutative residue of a projective pseudodifferential projection vanishes.
This also gives a very direct proof in the classical setting.
\end{abstract}



\maketitle

\section{Introduction}

The algebra of symbols of classical pseudodifferential operators $\PsDO_{cl}(X,E)$ on 
a closed manifold $X$ acting on sections of a vector bundle $E$ can be defined as the 
quotient of $\PsDO_{cl}(X,E)$ by the ideal of smoothing operators. Since pseudodifferential  
operators are smooth off the diagonal the symbol algebra is localized on the diagonal and 
it therefore can also be defined locally, using the product expansion formula and the change 
of charts formula for pseudodifferential operators. That the local heat kernel coefficients 
and the index of elliptic pseudodifferential operators are locally computable relies on the 
fact that the index and asymptotic spectral properties of pseudodifferential operators depend 
only on their class in the symbol algebra. Note that the principal symbol of a 
pseudodifferential operator is a section of the bundle of endomorphisms of $\pi^* E$, 
where $\pi: T^* X \to X$ is the canonical projection.

The bundle of endomorphisms of a complex hermitian vector bundle is a bundle of simple 
matrix algebras with $*$-structure. However, not all bundles of simple matrix algebras with 
$*$-structure, so-called Azumaya bundles, are isomorphic to endomorphism bundles of hermitian 
vector bundles. The obstruction is the so-called Dixmier-Douady class in $H^3(X,\mathbb{Z})$. 
Given an Azumaya bundle $\mathcal{A}$ it is possible to construct algebras of symbols whose 
principal symbols takes values in the space of sections of the Azumaya bundle 
$\pi^*\mathcal{A}$ (see for example \cite{MMS05} and the discussion in \cite{MMS06}). 
Following \cite{MMS05} we refer to such a symbol algebra as the algebra of symbols of 
projective pseudodifferential operators. For such symbol algebras one can define an index 
and Mathai, Melrose and Singer \cite{MMS06} proved an index formula for projective 
pseudodifferential operators, analogous to the Atiyah-Singer index formula. The topological 
index in this case is a map from twisted $K$-theory to $\mathbb{R}$. It has also been shown 
in \cite{MMS06} that any oriented manifold admits a projective Dirac operator even if the 
manifold does not admit a spin structure. In this case its index may fail to be an integer.

Another important quantity that depends only on the class of the symbol of a 
pseudodifferential operator is the so-called Wodzicki residue or noncommutative residue. 
Up to a factor it is the 
unique trace on the algebra of pseudodifferential operators. The Wodzicki residue appeared 
first as a residue of a zeta function measuring spectral asymmetry (\cite{APS76,Wo84}). 
Wodzicki showed that the regularity of the $\eta$-function of a Dirac-type operator at 
zero -- a necessary ingredient to define the $\eta$-invariant -- follows as a special case 
from the vanishing of the Wodzicki residue on pseudodifferential projections (as remarked 
by Br\"uning and Lesch \cite{BL99} the regularity of the $\eta$-function at zero for any 
Dirac type operator and the vanishing of the Wodzicki residue on pseudodifferential 
projections are actually equivalent). The regularity of the $\eta$-function was proved by 
Atiyah, Patodi and Singer in \cite{APS76} in the case when $X$ is odd dimensional and 
later by Gilkey (\cite{G81}) in the general case using $K$-theoretic arguments. Note that 
whereas the Wodzicki residue can be locally computed its vanishing on pseudodifferential 
projections is not a local phenomenon. Gilkey \cite{G79} constructed a pseudodifferential 
projection whose residue density is non-vanishing but integrates to zero. 

In our paper we show that the Wodzicki residue can also be defined for projective 
pseudodifferential operators (this has already been observed in \cite{MMS06}) and show 
that it vanishes on projections in case the dimension of the manifold is odd. 
Our proof is based on the  Leray-Hirsch theorem in twisted $K$-theory and a purely 
algebraic result on `residue-traces' in filtered rings. 

If $L$ is a filtered ring then we call a linear functional $\tau: L \to \mathbb{C}$ a 
residue trace if $\tau(L^{-N})=\{0\}$ for $N$ large enough. We prove that the value of 
$\tau$ on projections depends only on their class in $L^{(0)}:=L/L^{-1}$. Thus, if the 
map $K^0_{alg}(L) \to K^0_{alg}(L^{(0)})$ is surjective the map $\tau$ descends to a map 
from the algebraic $K$-theory of $L^{(0)}$ to $\mathbb{C}$. This result can be applied 
to the Wodzicki residue showing that it descends to a map from twisted $K$-theory  
$K^0(S^*X,\pi^*\mathcal{A})$ to $\mathbb{C}$. We then use the Leray-Hirsch theorem to 
show that this map actually vanishes. We reduce the 
problem to positive spectral projections of generalized Dirac operators for which it 
is known \cite{BG92} that the residue density vanishes.

\section{Convolution bundles and Azumaja bundles}

Pseudodifferential operators on a smooth closed Riemannian manifold $X$ acting on sections of a vector 
bundle $E$ can be understood as co-normal distributional sections in the vector bundle 
$E \boxtimes E^*$\footnote{$E \boxtimes E^*$ denotes the external tensor product of $E$ 
and its dual bundle $E^*$, i.e., the fibre over a point $(x,y)$ is $E_x\otimes E^*_y$.} 
over the space $X \times X$, by identifying the operators with their distributional kernel. 
The bundle $E \boxtimes E^*$ has the following structure that allows to convolve kernels of 
integral operators: any element in the fibre over ${(x,y)}$ may be multiplied by an element 
in the fibre over ${(y,z)}$ to give an element in the fibre over ${(x,z)}$. 
Moreover, this multiplication satisfies natural conditions such as associativity. 
In order to define projective pseudodifferential operators it is 
convenient to formalize this structure, as we shall do in this section.

\subsection{Convolution bundles}

Let $\mathcal{U}$ denote an open neighborhood of the diagonal $\Delta(X)$ 
in $X \times X$ which is symmetric under the reflection map $s:(x,y)\mapsto(y,x)$. 
Let $p_{ik}:X \times X \times X \to X \times X$ be defined by 
$p_{ik}(x_1,x_2,x_3)=(x_i,x_k)$ and set 
$\widetilde{\mathcal{U}}
   :=p_{12}^{-1}(\mathcal{U}) \cap p_{23}^{-1}(\mathcal{U})  
   \cap p_{13}^{-1}(\mathcal{U})$. 
Denote by $\widetilde p_{ik}$ the restriction of the map $p_{ik}$ to $\widetilde{\mathcal{U}}$.

\begin{definition}\label{def:bundle}
Let $\pi: F \to \mathcal{U}$ be a locally trivial vector bundle with typical fibre 
$\mathrm{Mat}(k)$, the complex $k\times k$-matrices. We call $F$ a \emph{convolution bundle} 
if there exists a homomorphism of vector bundles 
$m: \widetilde{p}_{12}^* F \otimes \widetilde{p}_{23}^* F \to F$ 
such that the following conditions are satisfied:
\begin{itemize}
 \item[(i)] The following diagram is commutative: 
   \begin{equation*}
    \xymatrix@1@=1.2cm@M=5pt{{\widetilde{p}_{12}^* F \otimes \widetilde{p}_{23}^* F} 
    \ar[d] \ar[r]^{\qquad m}  & F \ar[d] \\
    \widetilde{\mathcal{U}}   \ar[r]^{\widetilde{p}_{13}} & \mathcal{U}}
   \end{equation*}
 \item[(ii)] $m$ is associative, i.e., whenever $f_{ij}$ belong to the fibre 
   $F_{(x_i,x_j)}$ then 
   \begin{align*}
    m \left( m(f_{12}\otimes f_{23}) \otimes f_{34}\right)= 
    m \left( f_{12}\otimes m(f_{23} \otimes f_{34})\right).\footnotemark  
   \end{align*} 
 \item[(iii)] There\footnotetext{Here we implicitely assume that $(x_1,x_2,x_3)$, 
   $(x_1,x_3,x_4)$, and $(x_1,x_2,x_4)$ belong to $\widetilde{\mathcal{U}}$.} 
   is an atlas $\{\mathcal{O}_\alpha\}$ of $\mathcal{U}$ together with local trivializations 
    $$\phi_\alpha: \pi^{-1}\mathcal{O}_\alpha \to \mathcal{O}_\alpha \times 
      \mathrm{Mat}(k),$$ 
   such that
    $$\phi_\alpha(m(f_{12} \otimes f_{23})) = \phi_\alpha(f_{12}) \cdot \phi_\alpha(f_{23})$$
   whenever $f_{ij}\in F_{(x_i,x_j)}$ with 
   $(x_1,x_2,x_3)\in \widetilde{p}_{12}^{-1}(\mathcal{O}_\alpha)
   \cap\widetilde{p}_{23}^{-1}(\mathcal{O}_\alpha)$. 
 \end{itemize}
\end{definition}

\begin{definition}\label{def:star}
A \emph{$*$-structure} on $F$ is a conjugate linear map $*: F \to F$ of vector bundles 
such that 
\begin{gather*}
 \xymatrix@1@=1.2cm@M=5pt{
     F  \ar[d] \ar[r]^{*} & F \ar[d] \\
    \mathcal{U} \ar[r]^{s} & \mathcal{U} 
 }
\end{gather*}
commutes, such that
$
 (m(f \otimes g))^*=m(g^* \otimes f^*), 
$
and such that the above local trivializations additionally satisfy  
$$
 \forall f \in \pi^{-1}(\mathcal{O}_\alpha \cap s(\mathcal{O}_\alpha)): \qquad \phi_\alpha(f^*)=\phi_\alpha(f)^*, 
$$
where the star on the right hand side denotes the hermitian conjugation of matrices.
We will refer to a convolution bundle with $*$-structure as a \emph{$*$-convolution bundle}.
\end{definition}

Note that $E \boxtimes E^*$ is a particular example for a $*$-convolution bundle; in this 
case we can choose $\mathcal{U}=X \times X$. The restriction of a $*$-convolution 
bundle $F$ to the diagonal in $X \times X$ is a bundle $\mathcal{A}$ of finite dimensional 
simple $C^*$-algebras. Following the literature 
we refer to such bundles of matrix algebras as Azumaja bundles. 

As shown in \cite{MMS05} any Azumaja bundle $\mathcal{A}$ 
on $X$ gives rise to a convolution bundle near the diagonal in the following way, using an 
atlas of local trivializations with respect to a good 
cover\footnote{A cover is good if finite intersections of elements therein are either empty 
or contractible.} 
$\{{U}_\alpha\}$ of $X$: The transition functions $\sigma_{\alpha \beta}$ are smooth functions 
on ${U}_{\alpha \beta}={U}_{\alpha} \cap {U}_{\beta}$ with values in the automorphisms of 
$\mathrm{Mat}(k)$. Since all automorphims are inner we 
can choose local functions $\varphi_{\alpha \beta}: {U}_{\alpha \beta} \to SU(k)$
that implement $\sigma_{\alpha \beta}$, i.e., 
$\sigma_{\alpha \beta}(x)(A) = \varphi_{\alpha \beta}(x) A \varphi_{\alpha \beta}^{-1}(x)$. 
In general, the functions $\varphi_{\alpha \beta}$ may violate the co-cycle condition 
and therefore are not the transition functions of a vector bundle. The cocycle condition 
for the $\sigma_{\alpha\beta}$ together with the condition that the $\varphi_{\alpha \beta}$ 
are chosen in $SU(k)$ show that any 
$\varphi_{\alpha \beta} \varphi_{\beta \gamma} \varphi_{\gamma \alpha}$
must be a constant function on ${U}_\alpha \cap {U}_\beta \cap {U}_\gamma$, 
equal to an $k$-th root of unity times the identity 
matrix.\footnote{On different triple intersections, the resulting unit-root can be different. 
This induces a torsion element in $H^3(X,\mathbb{Z})$, the Dixmier-Douady class.} 
Then we obtain a convolution bundle $F$ with typical fibre $\mathrm{Mat}(k)$ 
on a neighborhood of the diagonal by choosing the transition functions
 $$\phi_{\alpha\beta}(x,y)(A)=\varphi_{\alpha \beta}(x) A\, \varphi_{\alpha \beta}(y)^{-1},
   \qquad A\in \mathrm{Mat}(k),$$ 
on ${U}_{\alpha\beta} \times {U}_{\alpha\beta}$. 
There are also other possible extensions of $\mathcal{A}$, cf.\ \cite{MMS06}, and Proposition 
\ref{prop:transition}, below. 

\begin{remark}\label{rem:atlas}
In the sequel it will be occasionally convenient to choose an atlas for $F$ consisting of sets 
$\mathcal{O}_\alpha:=U_\alpha\times U_\alpha$, where $\{U_\alpha\}$ is a good cover of $X$; 
the corresponding trivialisations we shall denote by $\phi_\alpha$ $($so we use the same notation 
as in Definition {\rm\ref{def:bundle}.(iii)} above, but possibly have changed the atlas$)$. 
\end{remark}
\subsection{Transition functions} 

In the previous section we have seen how an Azumaja bundle leads to a convolution bundle by 
choosing certain transition functions. Let us now have a closer look to the transition functions 
of an arbitrary $*$-convolution bundle. Fix an atlas as explained in Remark \ref{rem:atlas} and let 
$\phi_{\alpha\beta}:\mathcal{O}_\alpha\cap \mathcal{O}_\beta=:\mathcal{O}_{\alpha\beta}
\lra GL(\mathrm{Mat}(k))$
be the transition functions defined by 
\begin{equation*}\label{eq:trans0}
 \phi_\beta\circ\phi_\alpha^{-1}\big((x,y),A\big)=\big((x,y),\phi_{\alpha\beta}(x,y)(A)\big),
\end{equation*}
Then condition iii$)$ of Definition \ref{def:bundle} is equivalent to 
 $$\phi_{\alpha\beta}(x,y)(A)\phi_{\alpha\beta}(y,z)(B)=\phi_{\alpha\beta}(x,z)(AB).$$
In particular, 
\begin{equation}\label{eq:trans1}
 (x,x)\mapsto\phi_{\alpha\beta}(x,x):\mathcal{O}_{\alpha\beta}\cap\Delta(X)
 \lra\mathrm{Aut}(\mathrm{Mat}(k)).
\end{equation} 
Moreover, Definition \ref{def:star} on the level of the transition functions means that 
\begin{equation}\label{eq:trans2a}
 \phi_{\alpha\beta}(x,y)(A^*)=\phi_{\alpha\beta}(y,x)(A)^*.
\end{equation}

\begin{proposition}\label{prop:transition}
Let $F$ be a $*$-convolution bundle with transition functions $\phi_{\alpha\beta}$ as described above. 
Then
\begin{equation}\label{eq:trans2}
 \phi_{\alpha\beta}(x,y)(A)
 =\lambda_{\alpha\beta}(x,y)\varphi_{\alpha\beta}(x)A\,\varphi_{\alpha\beta}(y)^{-1}
\end{equation}
with mappings 
 $$\varphi_{\alpha\beta}:\mathcal{O}_{\alpha\beta}\lra SU(k), \qquad 
   \lambda_{\alpha\beta}:\mathcal{O}_{\alpha\beta}\lra \cz,$$
satifying  
 $$\lambda_{\alpha\beta}(x,x)=1,\qquad 
   \lambda_{\alpha\beta}(x,y)\lambda_{\alpha\beta}(y,z)=\lambda_{\alpha\beta}(x,z),\qquad
   \lambda_{\alpha\beta}(x,y)=\overline{\lambda_{\alpha\beta}(y,x)},$$
and such that all $\varphi_{\alpha\beta}\varphi_{\beta\gamma}\varphi_{\gamma\alpha}$  
are constant functions on their domain of definition, equal to a $k$-th root of unity times 
the identity matrix.
\end{proposition} 
\begin{proof}
Combining \eqref{eq:trans1} with \eqref{eq:trans2a} we find $\varphi_{\alpha\beta}$ with 
\begin{equation*}
 \phi_{\alpha\beta}(x,x)(A)=\varphi_{\alpha\beta}(x)A\,\varphi_{\alpha\beta}(x)^{-1}, 
\end{equation*}
since all automorphisms of $\mathrm{Mat}(k)$ are inner. 
Now let us define
 $$\phi_{\alpha\beta}^\prime(x,y)(A)
   =\varphi_{\alpha\beta}(x)^{-1}\phi_{\alpha\beta}(x,y)(A)\varphi_{\alpha\beta}(y).$$
We then have 
 $$\phi_{\alpha\beta}^\prime(x,x)(A)=A,\qquad 
   \phi_{\alpha\beta}^\prime(A)(x,y)\phi_{\alpha\beta}^\prime(y,z)(B)=\phi_{\alpha\beta}^\prime(x,z)(AB).$$
It follows that 
$\phi_{\alpha\beta}^\prime(x,y)(AB)=\phi_{\alpha\beta}^\prime(x,y)(A)\phi_{\alpha\beta}^\prime(y,y)(B)
=\phi_{\alpha\beta}^\prime(x,y)(A)B$
and, analogously, $\phi_{\alpha\beta}^\prime(x,y)(AB)=A\phi_{\alpha\beta}^\prime(x,y)(B)$. 
Therefore, for all matrices $A$, 
 $$\phi_{\alpha\beta}^\prime(x,y)(\mathbf{1})A=\phi_{\alpha\beta}^\prime(x,y)(A)
   =A\phi_{\alpha\beta}^\prime(x,y)(\mathbf{1}),$$
where $\mathbf{1}$ is the identity matrix. 
This shows $\phi_{\alpha\beta}^\prime(x,y)(\mathbf{1})$ is a multiple of the identity matrix. 
Denoting the corresponding factor by $\lambda_{\alpha\beta}(x,y)$, the claim follows. 
\end{proof}

\section{Projective Pseudodifferential Operators}

Projective pseudodifferential operators have been defined in \cite{MMS05}. We adapt this definition to fit 
in our setting of convolution bundles.

\subsection{Pseudodifferential operators}

To clarify notation let us briefly recall the definition of classical (or polyhomogeneous) 
pseudodifferential operators on an open subset $\Omega$ of $\rz^n$. Let $V\cong\cz^k$ be a 
$k$-dimensional vector space.

A symbol of order $m\in\rz$ is a smooth function 
$a:\Omega\times\Omega\times\rz^n\to\mathrm{End}(V)=V\otimes V^*$ satisfying estimates 
 $$\left\|\partial^\alpha_\xi\partial^\beta_{(x,y)} a(x,y,\xi)\right\|\le 
   C_{\alpha\beta K}(1+|\xi|)^{m-|\alpha|}$$
for any multi-indices $\alpha,\beta$ and any compact subset $K$ of $\Omega\times\Omega$, 
and having an asymptotic expansion $a\sim\sum\limits_{j=0}^\infty \chi a_{m-j}$ with a
zero-excision function $\chi=\chi(\xi)$ and homogeneous components $a_{m-j}$, i.e., 
 $$a_{m-j}(x,y,t\xi)=t^{m-j}a_{m-j}(x,y,\xi)$$
for all $(x,\xi)$ with $\xi\not=0$ and all $t>0$. The pseudodifferential operator 
$\op(a):C^\infty_{0}(\Omega,V)\to C^\infty(\Omega,V)$ associated with $a$ is 
 $$[\mathrm{op}(a)\varphi](x)=\iint e^{i(x-y)\xi}a(x,y,\xi)\varphi(y)\,dy\dbar\xi,
   \qquad \varphi\in C^\infty_0(\Omega,V).$$ 
An operator $R:C^\infty_{0}(\Omega,V)\to C^\infty(\Omega,V)$ is called smoothing if it has 
a smooth integral kernel $k\in C^\infty(\Omega\times\Omega,\mathrm{End}(V))$, i.e., 
 $$(R\varphi)(x)=\int_\Omega k(x,y)\varphi(y)\,dy,\qquad \varphi\in C^\infty_0(\Omega,V).$$
A pseudodifferential operator of order $m\in\rz$ on $\Omega$ is an operator of the form 
$A=\op(a)+R$, where $a$ is a symbol of order $m$ and $R$ is smoothing. 

Any pseudodifferential operator $A=\op(a)+R$ of order $m$ can be represented in the form 
$\op(a_L)+R^\prime$, where $a_L(x,\xi)$ is a $y$-independent `left-symbol' of order $m$; up to 
order $-\infty$ the left-symbol is uniquely determined by the asymptotic expansion 
 $$a_L(x,\xi)\sim \sum_{|\alpha|=0}^\infty\frac{1}{\alpha!}\partial^\alpha_\xi D^\alpha_y 
   a(x,y,\xi)\Big|_{x=y}.$$
The homogeneous components of $A$ are by definition those of $a_L$,
 $$\sigma_{m-j}(A)(x,\xi):=(a_L)_{m-j}(x,\xi).$$  
By the Schwarz kernel theorem, we can identify $A$ with its distributional 
kernel 
 $$K_A\in \scrD^\prime(\Omega\times\Omega,V\otimes V^*),$$
the topological dual of $C^\infty_0(\Omega\times\Omega;V^*\otimes V)$. 
It is uniquely defined by the relation 
 $$\skp{K_A}{\psi\otimes\varphi}=\skp{\psi}{A\varphi},\qquad 
   \psi\in C^\infty_0(\Omega,V^*),\quad
   \varphi\in C^\infty_0(\Omega,V).$$
Denoting by $\mathrm{tr}:V^*\otimes V\to\cz$ the canonical contraction map, we have explicitly 
 $$\skp{K_A}{u}=\int_\Omega \mathrm{tr}[A u(x,\cdot)](x)\,dx,\qquad 
   u\in C^\infty_0(\Omega\times\Omega;V^*\otimes V).$$
By pseudo-locality, 
$K_A\in C^\infty(\Omega\times\Omega\setminus\Delta(\Omega),V\otimes V^*)$. 

If $U\subset X$ is a coordinate neighborhood, we can pull-back the local operators 
under the coordinate map. The resulting space of operators we shall denote by 
$\PsDO^m_{cl}(U;\mathrm{End}(V))$, the subspace of smoothing operators by  
$\PsDO^{-\infty}(U;\mathrm{End}(V))$. 

\subsection{Projective pseudodifferential operators}

In the following choose an atlas as explained in Remark \ref{rem:atlas}. 

\begin{definition}\label{def:pseudo}
 Let $F$ be a $*$-convolution bundle over $\calU$. A distribution 
 $A \in \calD^\prime(\calU, F)$ 
 is called a projective pseudodifferential operator of order $m\in\rz$ if 
 \begin{itemize}
  \item[(i)] $A$ is smooth outside the diagonal, 
  \item[(ii)] for any $\alpha$ the distribution 
   $\left( \phi_\alpha^{-1} \right)^* A\big|_{U_\alpha\times U_\alpha}$ is the distributional kernel of a 
   pseudodifferential operator $A_\alpha\in\PsDO^m_{cl}(U_\alpha;\mathrm{End}(\cz^{k}))$. 
 \end{itemize}
 We denote the vector space of $m$-th order projective pseudodifferential 
 operators by $\PsDO^m_{cl}(\mathcal{U};F)$, the subspace of smoothing elements by 
 $\PsDO^{-\infty}(\mathcal{U};F)$.  
\end{definition}

The subspace $\mathrm{Diff}^m(\mathcal{U};F)$ of projective differential operators consists 
of all projective pseudodifferential operators which are supported on the diagonal. 

\begin{remark}
If $\calU = X \times X$ and $F=E\boxtimes E^*$ for a bundle $E$ over $X$ 
then $\PsDO^m_{cl}(\mathcal{U};F)$ coincides with $\PsDO^m_{cl}(X;E,E)$, the 
pseudodifferential operators of order $m$ acting on sections into $E$. 
\end{remark}

Though projective pseudodifferential operators, in general, are not operators in the 
usual sense (i.e., acting between sections of vector bundles) all elements of 
the standard calculus can be generalized to this setting. In particular, 
the $*$-structure gives rise to a conjugation on $\PsDO^m_{cl}(\mathcal{U};F)$, 
defined by $A^*(x,y):=(A(y,x))^*$ in the distributional sense. 

Let $A$ be a projective pseudodifferential operator with local representatives $A_\alpha$ and $A_\beta$, 
cf.\ Definition \ref{def:pseudo}, where $\mathcal{O}_{\alpha\beta}$ is not empty. 
By passing to local coordinates on $U_\alpha\cap U_\beta$, we can associate with $A_\alpha$ 
and $A_\beta$ local symbols $a_\alpha(x,\xi)$ and $a_\beta(x,\xi)$, respectively. These symbols are 
then related by 
\begin{align}\label{eq:asymptotic}
\begin{split}
 a_\beta(x,\xi)
 &=\sum_{|\gamma|=0}^\infty 
   \frac{1}{\gamma!}\partial^\gamma_\xi D^\gamma_y\Big|_{y=x}
   \phi_{\alpha\beta}(x,y)\big(a_\alpha(x,\xi)\big)\\ 
 &=\sum_{|\gamma|=0}^\infty 
   \frac{1}{\gamma!}\partial^\gamma_\xi D^\gamma_y\Big|_{y=x}
   \big[\lambda_{\alpha\beta}(x,y)\varphi_{\alpha\beta}(x)a_\alpha(x,\xi)\varphi_{\alpha\beta}(y)^{-1}\big],
\end{split}
\end{align}
where the transition function $\phi_{\alpha\beta}$ is as described in \eqref{eq:trans0} and 
Proposition \ref{prop:transition}. Note that this behaviour, in general, differs from the standard case, 
due the factor $\lambda_{\alpha\beta}(x,y)$. 
However, \eqref{eq:asymptotic} together with $\lambda_{\alpha\beta}(x,x)=1$ shows that with $A$ we can 
associate a well-defined homogeneous principal symbol 
 $$\sigma_m(A)(x,\xi)\in C^\infty(S^*X,\pi^*\mathcal{A}),$$
where $\pi:S^*X\to X$ is the canonical co-sphere bundle over $X$. Vice versa, any given such 
section can be realized as the principal symbol of a projective pseudodifferential operator. 

If the projective pseudodifferential operators $A_1$ and $A_2$ are supported in a 
sufficiently small neighborhood of the diagonal in $\calU$ their usual composition 
$$ 
 (A_1 \circ A_2)(x,z)= \int_{X} m\left(A_1(x,y) \otimes A_2(y,z)\right)\, dy
$$
is a distribution. By passing to local coordinates and using the composition theorems for 
pseudodifferential operators one can see that $A_1 \circ A_2$ is a projective pseudodifferential 
operator. The homogeneous principal symbol behaves multiplicative under composition. Of course, 
any projective pseudodifferential operator can be written as a sum of two operators, where one 
is smoothing and the other is supported near the diagonal. Summarizing, the coset space 
\begin{equation}\label{eq:coset}
 L^*_{cl}(\calU,F):=\PsDO^*_{cl}(\calU,F) / \PsDO^{-\infty}(\calU,F)
\end{equation}
is a filtered $*$-algebra. As in the standard case, asymptotic summations of sequences 
of projective operators of one-step decreasing orders are possible and parametrices 
(i.e., inverses modulo smoothing remainders) to elliptic elements can be constructed.  

\begin{theorem}\label{wres}
 Let $F$ be a $*$-convolution bundle and let $A$ be a projective pseudodifferential 
 operator. For $x\in X$ define
 $$
  \mathrm{WRes}_x(A) := \int_{S^*_xX} \tr a_{-n}(x,\xi) \,d\sigma(\xi)\,dx,
 $$ 
 where $a_{-n}(x,\xi)$, $n=\mathrm{dim}\,X$, is the homogeneous component of 
 order $-n$ of a symbol of a local representative $A_\alpha$ with $x\in\mathcal{O}_\alpha$, 
 cf. Definition {\rm\ref{def:pseudo}}. Then  $\mathrm{WRes}_x(A)$ is well-defined and 
 defines a global density on $X$. Moreover,
 $$ 
  \mathrm{WRes}(A) := \int_X \mathrm{WRes}_x(A) 
 $$
 defines a trace functional on the algebra $L^*_{cl}(\calU,F)$, the so-called 
 \emph{noncommutative residue} or Wodzicki residue.
\end{theorem}
\begin{proof}
 Let $A_\beta$ be another local representative and $x\in\mathcal{O}_\beta$. Fixing local coordinates on 
 $\mathcal{O}_\alpha\cap\mathcal{O}_\beta$, the local symbols $a_\alpha$ and $a_\beta$ are related by 
 the asymptotic expansion \eqref{eq:asymptotic}. Following the proof in \cite{FGLS} terms containing 
 a derivative $\partial_\xi^\gamma$, $|\gamma|\ge1$, vanish under integration. We thus obtain the same 
 value for $\mathrm{WRes}_x(A)$ using either $a_\alpha(x,\xi)$ or $a_\beta(x,\xi)$. 
 That $\mathrm{WRes}_x(A)$ transforms as density under changes of 
 coordinates is seen as in the standard case, cf.\ \cite{FGLS}. 
 
 To see that the integral of the residue density defines a  trace functional we need to show that it 
 vanishes on commutators $[A,B]$. To this end fix a cover $\{U_\sigma^\prime\}$ of $X$ by coordinate maps 
 together with a subordinate partition of unity, such that $U_\sigma^\prime\cup U_\rho^\prime$ is contained 
 in some $U_\alpha$ whenever $U_\sigma^\prime\cap U_\rho^\prime$ is not empty. We then 
 can write $A= \sum_\sigma A_{\sigma}$ and $B= \sum_\sigma B_{\sigma}$ modulo smoothing operators, 
 where the $A_\sigma$ and $B_\sigma$ are supported in 
 $\mathcal{O}_{\sigma}^\prime:=U_\sigma^\prime \times U_\sigma^\prime$. Then the commutator
 $[A,B]$ can be written as a sum of terms $[A_\sigma,B_\rho]$. Such a commutator is smoothing 
 if $\mathcal{O}_{\sigma}^\prime\cap\mathcal{O}_{\rho}^\prime$ is empty. Otherwise it 
 is contained in some set $\mathcal{O}_\alpha$. Therefore the calculation reduces to a local one, which 
 is not different from the one for usual pseudodifferential operators that can be found in \cite{FGLS}. 
\end{proof}

For purposes below let us establish the following result: 

\begin{proposition}\label{projection}
Let $F$ be a $*$-convolution bundle and $\mathcal{A}$ be the Azumaja 
bundle obtained by restricting $F$ to the diagonal. Moreover, let 
$p\in C^\infty(S^*X,\pi^*\mathcal{A})$ with $p^2=p$. Then there exists a 
projective pseudodifferential operator $P\in L^0_{cl}(\calU,F)$ which is a projection, 
i.e., $P^2=P$, and which has $p$ as its principal symbol. If, additionally, $p^*=p$ 
then $P$ can be chosen such that $P^*=P$. 
\end{proposition}
\begin{proof}
If $p^2=p$ then $e=2p-1$ is an idempotent. We now construct a projective 
pseudodifferential operator $E\in L^0_{cl}(\calU,F)$ which is an idempotent and has 
$e$ as its principal symbol. 
Then $P=(1+E)/2$ is the desired projection. Let $\wt{E}\in L^0_{cl}(\calU,F)$ be any 
element having $e$ as principal symbol. Then $\wt{E}^2=1-R$ with a remainder 
$R\in L^{-1}_{cl}(\calU,F)$. If $\sum_{k=0}^\infty c_kr^k$ denotes the Taylor series of 
$f(r)=1/\sqrt{1-r}$ let $S\in L^0_{cl}(\calU,F)$ have the asymptotic expansion 
$\sum\limits_{k=0}^\infty c_kR^k$. Then $(1-R)S^2=1$ and $S$ commutes with 
$\wt{E}$, since $R$ does. Then define $E=\wt{E}S$. 

In case  also $p^*=p$, first choose $E_0$ having $e$ as principal symbol. Then set 
$\wt{E}=E_0E_0^*$ and proceed as before; note that $R^*=R$ and hence $S^*=S$. 
\forget{
The proof follows \cite{Schu01}. Using local coordinates together with a partition of 
unity we can construct a $P_1\in L^0_{cl}(\calU,F)$ having $p$ as principal symbol. 
If $V:=\cz\setminus\{0,1\}$ then $\lambda-P_1$ is elliptic for any 
$\lambda\in V$. Thus there is a parametrix $Q(\lambda)$, depending holomorphically 
on the parameter $\lambda\in V$ (i.e., the local symbols depend holomorphically on $\lambda$). 
Then one takes   
 $$P=\int_\scrC Q(\lambda)\,d\lambda,$$
where $\scrC=\partial U_{1/2}(1)$ is the counter-clockwise oriented boundary of the disc of 
radius $1/2$ centred in $1$ $($more precisely, we decompose $Q(\lambda)$ in a sum of 
local terms and integrate each of these terms seperately over $\scrC)$. 
If also $p^*=p$ we repeat the above construction, replacing $P_1$ by $P^*P$. 
}
\end{proof}

\section{The noncommutative residue in twisted $K$-theory}

\subsection{Twisted K-theory}

Suppose that $\mathcal{A}$ is an Azumaja bundle over a compact manifold $X$.
The twisted K-theory is defined to be the K-theory of the $C^*$-algebra
of continuous sections $C(X; \mathcal{A})$ of $\mathcal{A}$.

If $Y \subset X$ is a closed subset then the set of sections $C(X,Y;\mathcal{A})$
vanishing on $Y$ is a closed two-sided ideal in $C(X,Y;\mathcal{A})$
and the quotient by this ideal can be identified with the space of continuous
sections $C(Y;\mathcal{A})$ of the Azumaja bundle $\mathcal{A}|_Y$. We therefore have
the six term exact sequence as a consequence of the six term exact sequence in the theory
of $C^*$-algebras.
\begin{gather*}
  \xymatrix@1@=1.2cm@M=2pt{
     K^0(X,Y;\calA) \ar[r] & K^0(X;\calA) \ar[r] & K^0(Y;\calA) \ar[d] \\
    K^1(Y;\calA)  \ar[u] \ar[r] & K^1(X;\calA) \ar[r] & K^1(X,Y;\calA) }
\end{gather*}
where the relative $K$-groups $K^*(X,Y;\calA)$ are defined as
$K_*(C(X,Y;\calA))$.

There is a natural map
\begin{gather*}
 K_*(C(X,Y;\calA)) \otimes_\mathbb{Z} K_*(C(X)) \mapsto K_*(C(X,Y;\calA) \hat \otimes C(X)).
\end{gather*}
Here $\hat \otimes$ is the tensor product of $C^*$-algebras which is well defined in this case as $C(X)$
is nuclear.
The usual multiplication
\begin{gather*}
 C(X,Y;\calA) \hat \otimes C(X) \to C(X,Y;\calA)
\end{gather*}
induces a map $K_*(C(X,\calA) \hat \otimes C(X)) \to K_*(C(X,\calA))$.
The composition of these two maps makes $K^*(X,Y;\calA)$
a module over the $\mathbb{Z}_2$-graded ring $K^*(X)$.
Choosing $Y=\emptyset$ defines a $K^*(X)$ module structure
on $K^*(X;\calA)$.  Note that the morphisms in the six term exact sequence are module
homomorphisms.

These observations can be used to prove the following Leray-Hirsch theorem: 

\begin{theorem}\label{leray}
 Let $R$ be a commutative torsion-free ring.
 Suppose that $\pi: M \xrightarrow{F} X$ is a compact smooth fibre bundle with fibre $F$ 
 over $X$ and let $\mathcal{A}$ be an Azumaja bundle over $X$. Assume that $K^*(F) \otimes_{\gz} R$ is a 
 free $R$-module and suppose there exist elements
 $c_1,\ldots,c_N \in K^*(M) \otimes_{\gz} R$ such that the $c_j|_{M_x}$ form a basis 
 for $K^*(M_x)\otimes_{\gz} R$ for every $x \in X$. Then the following map is an isomorphism$:$
 \begin{equation*}\label{eq:leray}
   K^*(X;\calA)\otimes_{\gz}  R^N \longrightarrow K^*(M,\pi^*(\calA)) \otimes_{\gz} R,\qquad 
   (p,\alpha) \mapsto \sum_{j=1}^N \alpha_j \pi^*(p) \cdot c_j. 
 \end{equation*}
\end{theorem}

Indeed, the usual proof  of the Leray-Hirsch theorem in topological $K$-theory (see e.g. \cite{hatcher}, Theorem 2.25)
can be adapted to our setting in the following way.
If $Y \subset X$ is a closed subset of $X$, we have the following diagram:
\begin{gather*}
  \xymatrix@=0.5cm@1@=0.5cm@=0.5cm@M=2pt{
     \ar[r] & K^*(X,Y;\calA) \otimes_{\gz} R^N \ar[d]^\Phi  \ar[r] &  K^*(X;\calA) \otimes_{\gz} R^N  
     \ar[d]^\Phi \ar[r] & K^*(Y;\calA)\otimes_{\gz} R^N   \ar[d]^\Phi  \ar[r] &  \\
     \ar[r] &K^*(\pi^{-1}X,\pi^{-1}Y;\calA) \otimes_{\gz} R \ar[r] & K^*(\pi^{-1}X;\calA)\otimes_{\gz} R \ar[r] & K^*(\pi^{-1}Y;\calA) \otimes_{\gz} R\ar[r] & }
\end{gather*}
Here  $\Phi$ is defined as in the theorem,
$\Phi(p,\alpha) = \sum\limits_{j=1}^N \alpha_j \pi^*(p) \cdot c_j$. The rows of this diagram are exact
since tensoring with $R^N$ and $R$ is an exact functor. All maps in the six term exact sequence are natural
and therefore the pull back $\pi^*$ commutes with them. Moreover, the maps in the six term exact sequence for
the pair $(\pi^{-1}X,\pi^{-1}Y)$ are $K^*(M)$ module homomorphisms. Thus, the diagram commutes.
Since $X$ is a finite cell complex one can proceed in the usual way using the $5$-lemma and induction
in the number of cells and the dimension  to prove the theorem.

\subsection{Residue-traces on filtered rings}

Let $L$ be a ring with filtration, i.e., $L=L^0\supset L^{-1}\supset L^{-2}\supset\ldots$ 
with sub-rings $L^{-j}$ and the multiplication induces maps $L^{-i}\times L^{-j}\to L^{-i-j}$ 
for any choice of $i,j$. 

A trace functional on $L$ is a map $\tau:L\to V$ for some vector space $V$ having the following properties: 
\begin{itemize}
 \item[$(1)$] $\tau$ is linear, $\tau(A+B)=\tau(A)+\tau(B)$ for all $A,B\in L$, 
 \item[$(2)$] $\tau$ vanishes on commutators, $\tau([A,B])=\tau(AB-BA)=0$ for all $A,B\in L$.
\end{itemize}
We call $\tau$ a residue-trace if, additionally,    
\begin{itemize}
 \item[$(3)$] there exists an $N$ such that $\tau(A)=0$ for all $A\in L^{-N}$. 
\end{itemize}
We shall now show that a residue-trace restricted to the set of projections in $L$ is insensible 
for lower order terms. The proof is elementary and purely algebraic. 

\begin{theorem}\label{prop1}
Let $\tau$ be a residue-trace on $L$ and $P,\wt{P}\in L$ be two projections, i.e., 
$P^2=P$ and $\wt{P}^2=\wt{P}$. If $P-\wt{P}\in L^{-1}$ then $\tau(P)=\tau(\wt{P})$. 
\end{theorem}
\begin{proof}
Set $R=\wt{P}-P$ and then define  
 $$A=P R P,\quad B=P R(1-P),\quad C=(1-P)R P,\quad D=(1-P)R(1-P).$$
Obviously then $\wt{P}=P+R$ and $R=A+B+C+D$. Using that $P(1-P)=(1-P)P=0$ we obtain 
\begin{align*}
 (P+R)(P+R)=&P+2A+B+C+A^2+AB+BC+ \\
 &+BD+CA+CB+DC+D^2.
\end{align*}
On the other hand, using that $\wt{P}$ is a projection, 
 $$(P+R)(P+R)=(P+R)=P+A+B+C+D.$$
Equating these two expressions and rearranging of terms yields 
 $$A^2+A+BC+D^2-D+CB+AB+BD+CA+DC= 0.$$
Multiplying this identity from the left and the right with $P$ and $1-P$, respectively, 
yields 
\begin{equation*}
 A^2+A+BC = 0,\qquad D^2-D+CB =0.
\end{equation*}
The first identity shows $A\in L^{-2}$ and $A=-BC$ modulo $L^{-4}$. 
Let us now rewrite these equations as 
\begin{equation*}
 A(1+A)=-BC,\qquad (-D)\big(1+(-D)\big)=-CB.  
\end{equation*}
Multiplying the first equation by $(1-A)$ yields $A\equiv -BC-(BC)^2$ modulo $L^{-6}$. Multiplying it 
with $(1-A+A^2)$ then yields $A\equiv-BC-(BC)^2-2(BC)^3$ modulo $L^{-8}$. Proceeding by induction we obtain 
 $$A\equiv \sum_{k=1}^\ell c_{k\ell}(BC)^{k}\mod L^{-2(\ell+1)}$$
for any $\ell\in\nz$ with suitable constants $c_{k\ell}$. In the same way, with the same 
constants $c_{k\ell}$,  
 $$-D\equiv \sum_{k=1}^\ell c_{k\ell}(CB)^{k}\mod L^{-2(\ell+1)}.$$
Therefore we have 
 $$A+D=\sum_{k=1}^\ell c_{k\ell}\left[B,(CB)^{k-1}C\right]\mod L^{-2(\ell+1)}.$$
Choosing $\ell$ large enough, we deduce that $\tau(A)+\tau(D)=0$. Furthermore, 
 $$\tau(B)=\tau(P R(1-P))=\tau((1-P)P R)=0$$ 
and, analogously, $\res(C)=0$. Altogether we obtain 
 $$\tau(\wt{P})=\tau(P)+\tau(A)+\tau(D)+\tau(B)+\tau(C)=\tau(P)$$
which is the claim we wanted to prove. 
\end{proof}

\subsection{The noncommutative residue}

We shall show that the noncommutative residue induces a map on twisted $K$-theory. 

\begin{proposition}
Let $\mathcal{A}$ be the Azumaja bundle obtained by restricting a $*$-con\-vo\-lution bundle 
$F$ to the diagonal. The noncommutative residue from Theorem {\rm\ref{wres}} descends to a 
group homomorphism
\begin{equation}\label{eq:wres}
 \mathrm{WRes}:K^0(S^*X,\pi^*\mathcal{A}) \to \mathbb{C},
\end{equation}
where $\pi:S^*X\to X$ denotes the co-sphere bundle over $X$. 
\end{proposition}
\begin{proof}
A typical element in $K^0(S^*X,\pi^*\mathcal{A})$ can be represented by a section  
$p\in C^\infty(S^*X,\pi^*\mathrm{Mat}_N(\mathcal{A}))$ which is $($pointwise$)$ a projection. 
This is possible, since the natural inclusion of the $K$-theory of the local $C^*$-algebra 
$C^\infty(S^*X,\pi^*\mathcal{A})$ into that of $C(S^*X,\pi^*\mathcal{A})$ is an isomorphism, 
cf. \cite{Bl98}. By Proposition \ref{projection} each such section is the principal symbol 
of a projective 
pseudodifferential operator $P\in L^0(X;\mathrm{Mat}_N(F))$ which is a projection. 
The noncommutative residue of the $K$-group element is then defined as $\mathrm{WRes}(P)$ 
in the sense of Theorem \ref{wres}. 

We have to show that this map is well-defined. So let 
$\wt{p}\in C^\infty(S^*X,\pi^*\mathrm{Mat}_M(\mathcal{A}))$ 
represent the same element as $p$ does. Let $\wt{P}\in L^0(X;\mathrm{Mat}_M(F))$ 
be associated with $\wt{p}$. 
Since $p$ and $\wt{p}$ are equivalent there exists a unitary 
$u\in C^\infty(S^*X,\pi^*\mathrm{Mat}_{M+N}(\mathcal{A}))$ such that 
$u(p\oplus 0_{M})u^{-1}$ coincides with $0_{N}\oplus\wt{p}$. 
Let $U\in L^0(X;\mathrm{Mat}_{M+N}(F))$ have $u$ as its principal symbol. Then 
 $$\mathrm{WRes}(U(P\oplus 0_{M})U^{-1})=\mathrm{WRes}(P\oplus 0_{M})=\mathrm{WRes}(P).$$
On the other hand $U(P\oplus 0_{M})U^{-1}$ is a projection having the same principal 
symbol as $0_{N}\oplus\wt{P}$. Thus, by Theorem \ref{prop1},  
 $$\mathrm{WRes}(U(P\oplus 0_{M})U^{-1}) =\mathrm{WRes}(0_{N}\oplus\wt{P})=\mathrm{WRes}(\wt{P}).$$
This shows that the noncommutative residue is independent of the choice of the representative. 
\end{proof}

\section{Twisted Dirac operators and connections}

Let $\mathcal{A}$ be the Azumaja bundle obtained by restricting a $*$-convolution bundle 
$F$ to the diagonal. 

\begin{definition}\label{def:connection}
A projective connection $\nabla=\nabla^F$ on $F$ is a linear map 
 $$Y\mapsto\nabla_Y:\quad C^\infty(X;TX)\lra \mathrm{Diff}^1(\mathcal{U};F)$$
satisfying, for any vector field $Y\in C^\infty(X,TX)$ and any function 
$f\in C^\infty(X)$,  
\begin{itemize}
 \item[$(1)$] $\nabla_{fY}=f\nabla_Y$, 
 \item[$(2)$] $[ \nabla_Y , f] = Yf$ for any $f\in C^\infty(X)$.
\end{itemize}
It is called a hermitian connection if additionally 
\begin{itemize}
 \item[$(3)$] $\nabla_Y^* + \nabla_Y + \mathrm{div}\, Y =0$
\end{itemize}
$($here, $f$ and $\mathrm{div}\,Y$ are considered as elements of $\mathrm{Diff}^0(\mathcal{U};F))$. 
\end{definition}

Note that in case $\mathcal{U}=X\times X$ and $F=E \boxtimes E^*$ for a vector bundle 
$E$ over $X$ we just recover a usual hermitian connection on $E$. One can always construct 
a projective hermitian connection from local hermitian connections 
by gluing with a partition of unity.

If $\nabla=\nabla^F$ is a projective connection and $\phi_\alpha$ is a local trivialization
of $F$ over $U_\alpha\times U_\alpha$ as described in Remark \ref{rem:atlas}, the corresponding 
local differential operator 
 $$\nabla^\alpha_Y\in\mathrm{Diff}^m(U_\alpha,\mathrm{End}(\cz^k))$$ 
is of the form  
 $$\nabla^\alpha_Y=Y+\Gamma_Y^\alpha(x),\qquad 
   \Gamma_Y^\alpha\in C^\infty(U_\alpha,\mathrm{Mat}(k)).$$ 
If we use another trivialisation $\phi_\beta$ of $F$ on $U_\beta\times U_\beta$, 
we have the relation 
 $$\Gamma_Y^\beta(x)=\phi_{\alpha\beta}(x,x)(\Gamma_Y^\alpha(x))
   +Y_y\phi_{\alpha\beta}(x,y)(\mathbf{1})\big|_{y=x},\qquad x\in U_\alpha\cap U_\beta,$$
where $\mathbf{1}$ is the identity matrix. Thus, analogous to the theory of standard connections, 
we may describe projective connections by `connection matrices' $\Gamma_Y^\alpha$ associated to 
a covering $X=\mathop{\mbox{\large$\cup$}}_\alpha U_\alpha$ satisfying the above compatibility 
relations. For a hermitian connection the connection matrices also have to be skew-symmetric, 
$\Gamma_Y^\alpha(x)^*=-\Gamma_Y^\alpha(x)$.

Suppose now $S$ is a Clifford module over $X$ and let $\gamma$ denote the Clifford 
multiplication. Moreover, let $\nabla^{S}$ be a connection on $S$ which is compatible 
with the Clifford structure. Writing $\wt{F}:=S\boxtimes S^*$ it is easy to see that 
$F\otimes\wt{F}$ is a $*$-convolution bundle over $\mathcal{U}$ and we can define 
the projective hermitian connection 
 $$\nabla:=\nabla^F\otimes 1+1\otimes\nabla^{S}$$ 
by choosing the corresponding connection matrices as  
 $$\Gamma_Y^{F,\alpha}(x)\otimes 1+1\otimes \Gamma_Y^{S,\alpha}(x),\qquad x\in U_\alpha,$$
where the $U_\alpha$ are chosen in such a way that both $F$ and $\wt{F}$ are locally trivial 
over $U_\alpha\times U_\alpha$. 
Then we can define the twisted Dirac operator 
 $$D:=(1\otimes\gamma)\circ \nabla\in \mathrm{Diff}^1(\mathcal{U};F\otimes\wt{F});$$
in fact, in each local trivialisation $\nabla$ is a usual hermitean connection and we can compose 
it locally with $1\otimes\gamma$. 

\section{Vanishing of the Wodzicki residue}

\begin{theorem}\label{thm:vanish}
 If $X$ is an odd dimensional oriented manifold 
 the map $\mathrm{WRes}$ of \eqref{eq:wres} vanishes identically.
\end{theorem}

As a direct consequence, $\mathrm{WRes}(P)=0$ for any projection 
$P\in L^0(\mathcal{U},F)$. 

\begin{proof}[Proof of Theorem \ref{thm:vanish}]
Suppose the dimension of $X$ is $n=2 \ell -1$.
Let $S=\mathop{\mbox{$\oplus$}}\limits_{k\text{ even}}\Lambda^k(T^*X)$ denote the bundle of 
even-degree forms over $X$ and let $*: \Lambda^k(T^*X) \to \Lambda^{n-k}(T^*X)$ the Hodge star
operator and denote by $d$ and $\delta$ the exterior differential and the co-differential respectively.
Define the operator $D^S$ acting on sections of $S$ as
 $$D^S=i^{\ell} * (\delta+(-1)^{k+1}d)\quad\text{on $k$-forms}.$$
Then by Proposition 1.22 and 2.8 in \cite{BGV04} this is a generalized Dirac operator, 
where the Clifford action on $S$ is given by 
$\gamma(\xi)=i^{\ell} * \left(\mathrm{int}(\xi)+(-1)^{k+1}\mathrm{ext}(\xi)\right)$ for $\xi\in T^*X$
and the compatible connection is the Levi-Civita connection.
The principal symbol of $D^S$ restricted to the co-sphere bundle is a self-adjoint involution and the
projection $\sigma_+(D^S)=\frac{1}{2}(\sigma(D^S)+1)$ onto its $+1$ eigenspace defines an element
in $K^0(S^*X)$. It is well known (see for instance \cite{APS76}) that restriction of this element to each
co-sphere $S^*_x X$ equals $2^\ell$ times the Bott element on $S^{n-1}$ which together with the class of the trivial
line bundle generates $K^0(S^{n-1})$.

For notational convenience denote by $K_{\rz}^*(X)$ the groups $K^*(X) \otimes_{\gz} \rz$.
By Theorem \ref{leray} applied to the co-sphere bundle of $X$, any element of 
$K^0_{\rz}(S^*X,\pi^*\mathcal{A})$ can be represented in the form 
 $$\alpha_0\pi^*(p)\cdot [\mathbf{1}]+\alpha_1\pi^*(p)\cdot [\sigma_+(D^S)]$$
for some $\alpha_0,\alpha_1 \in \rz$ and some $p\in K^0(X,\mathcal{A})$. Here both
the class $[\mathbf{1}]$ of the trivial line bundle and the class $[\sigma_+(D^S)]$ are understood
as elements in $K^0_\rz(S^*X)$.
The elements in $\alpha_0\pi^*(p)\cdot [\mathbf{1}]$ can be represented by projections in 
$C^\infty(X;\mathrm{Mat}_N(\mathcal{A}))$.
Therefore, the noncommutative residue of these elements vanishes. 
It remains to show that this is also true for the second summand. 

To this end let $p$ be a projection in $\mathrm{Mat}_N(C^\infty(X;\mathcal{A}))$. 
Let us define the new convolution bundle $F_p$ 
having fibre $p(x)\mathrm{Mat}_N(F)_{(x,y)}p(y)\subset \mathrm{Mat}_N(F)_{(x,y)}$ in $(x,y)$. 
We now apply the above construction and build a twisted Dirac operator $D_p$ with respect to 
$F_p\otimes\wt{F}$, $\wt{F}=S\boxtimes S^*$. Then $\sigma_+(D_p)$ represents the class 
$\pi^*([p])\cdot[\sigma_+( D^S)]$ in $K^0(S^*X,\pi^*\mathcal{A})$. 

The projective differential operator $D_p$ can now be used to construct a certain 
projection $Q \in L^0_{cl}(X,F_p)$  which has principal symbol 
as $\sigma_+(D_p)$ on $S^* X$.
In the case of a Dirac type operator $D$ acting on a vector bundle the projection would just be the
operator $\frac{1}{2}(|D|^{-1} D+1)$. The symbol of this projection can be constructed from the a parametrix
of $D$ and this construction is local modulo smoothing operators. That is the full symbol of $\frac{1}{2}(|D|^{-1} D+1)$
modulo smoothing terms in local coordinates depends only on the full symbol of $D$ in these local coordinates.
Thus, the construction can be repeated for the operator $D_p$ to yield an element in $L^0_{cl}(X,F_p)$
which we denote by $Q$ or formally $\frac{1}{2}(|D_p|^{-1} D_p+1)$. By construction $[\sigma(Q)] \in K^0_\rz(S^*X;\mathcal{A})$ 
is equal to  $\pi^*([p]) \cdot [\sigma_+(D^S)]$.

In  \cite{BG92} (Theorem 3.4) Branson and Gilkey have used invariant theory to show that the residue \emph{density} 
of the positive spectral projection for any generalized Dirac operator vanishes identically. 
Locally, $D_p$ is a generalized Dirac operator
and since the construction of the residue density is local the residue density of  $Q$
vanishes as well.
So we can conclude that the noncommutative residue of  $Q$ vanishes which
completes our proof. 
\end{proof}

\noindent{\bf Acknowledgements.} The authors would like to thank
Thomas Schick for comments and for pointing out a gap in an earlier version
of this paper.


\bibliographystyle{amsalpha}

\begin{thebibliography}{APS75}

\bibitem[APS76]{APS76}
M.F. Atiyah, V.K. Patodi, and I.M. Singer, \emph{Spectral asymmetry and
  {R}iemannian geometry {\rm III}}, Math. Proc. Cambridge Philos. Soc. \textbf{79}
  (1976), no.1, 71--99. 

\bibitem[BGV04]{BGV04}
N. Berline, E. Getzler and M. Vergne, \emph{Heat kernels and {D}irac operators},
Grundlehren Text Editions, Springer-Verlag, Berlin, 2004.

\bibitem[Bl98]{Bl98}
B. Blackadar, \emph{{$K$}-theory for operator algebras},
Mathematical Sciences Research Institute Publications,
Vol. 5, Second. Ed., Cambridge University Press, 1998.

\bibitem[BG92]{BG92}
T.P. Branson and P.B. Gilkey, \emph{Residues of the eta function for an operator of {D}irac type},
J. Funct. Anal. \textbf{108} (1992), no. 1, 47--87.

\bibitem[BL99]{BL99}
J. Br{\"u}ning and M. Lesch, 
\emph{On the {$\eta$}-invariant of certain nonlocal boundary value problems}, 
Duke Math. J. \textbf{96} (1999), no. 2, 425--468. 

\bibitem[FGLS96]{FGLS} 
B.\ Fedosov, F.\ Golse, E. Leichtnam and E. Schrohe, 
\emph{The noncommutative residue for manifolds with boundary}.
J. Funct. Anal. \textbf{142} (1996), no. 1, 1--31.
    
\bibitem[Gi79]{G79}    
P. Gilkey, \emph{The residue of the local eta function at the origin}, 
Math. Ann. \textbf{240} (1979), no. 2, 183--189.

\bibitem[Gi81]{G81}    
P. Gilkey, \emph{The residue of the global $\eta$ function at the origin},  
Adv. in Math. \textbf{40} (1981), no. 3, 290--307.

\bibitem[H09]{hatcher}
A. Hatcher, \emph{Vector Bundles \& K-Theory}, unpublished manuscript, available
at {\tt http://www.math.cornell.edu/~hatcher/VBKT/VBpage.html}, version 2.1, May 2009.
 
 \bibitem[MMS05]{MMS05}
V. Mathai, R.B. Melrose and I.M. Singer, 
\emph{The index of projective families of elliptic operators},
Geom. Topol. \textbf{9} (2005), 341--373.
  
\bibitem[MMS06]{MMS06}
V. Mathai, R.B. Melrose and I.M. Singer, \emph{Fractional analytic index},
J. Differential Geom. \textbf{74} (2006), no. 2, 265--292.

\forget{
\bibitem[Sch01]{Schu01}
B.-W.\ Schulze, \emph{An algebra of boundary value problems not requiring 
Shapiro-Lopatinskij conditions}, J. Funct. Anal. \textbf{179} (2001), no. 2, 374--408.
}

\bibitem[Wo84]{Wo84}
M. Wodzicki, \emph{Local invariants of spectral asymmetry},
Invent. Math. \textbf{75} (1984), no. 1,143--177. 

\end{thebibliography}

\end{document}